\documentclass[12pt,draftcls,onecolumn]{IEEEtran}
\IEEEoverridecommandlockouts

\usepackage{amsmath}
\usepackage{graphicx}
\usepackage{amssymb}
\usepackage{color}
\usepackage{fullpage}
\usepackage{graphicx}
\usepackage{amssymb}
\usepackage{color}
\usepackage{authblk}
\usepackage{psfrag}

\newtheorem{theorem}{Theorem}

\newtheorem{corollary}{Corollary}

\begin{document}

\title{Decentralized Event-Triggered Consensus over Unreliable Communication Networks}

\author{Eloy Garcia\thanks{Corresponding Author: \ttfamily{elgarcia@infoscitex.com}. \\ This work has been supported in part by AFOSR LRIR No. 12RB07COR. \\ Eloy Garcia is a contractor (Infoscitex Corp.) with the Control Science Center of Excellence, Air Force Research Laboratory, Wright-Patterson AFB, OH 45433. \\ Yongcan Cao is with the Department of Electrical and Computer Engineering, University of Texas at San Antonio, San Antonio, TX 78249 \\ David Casbeer is with the Control Science Center of Excellence, Air Force Research Laboratory, Wright-Patterson AFB, OH 45433. \\ A shorter version of this document was submitted to Automatica.}, Yongcan Cao, and David W. Casbeer}

                                         
\maketitle 
\begin{abstract}
This article studies distributed event-triggered consensus over unreliable communication channels. Communication is unreliable in the sense that the broadcast channel from one agent to its neighbors can drop the event-triggered packets of information, where the transmitting agent is unaware that the packet was not received and the receiving agents have no knowledge of the transmitted packet. Additionally, packets that successfully arrive at their destination may suffer from time-varying communication delays. In this paper, we consider directed graphs, and we also relax the consistency on the packet dropouts and the delays. Relaxing consistency means that the delays and dropouts for a packet broadcast by one agent can be different for each receiving node. We show that even under this challenging scenario, agents can reach consensus asymptotically while reducing transmissions of measurements based on the proposed event-triggered consensus protocol. In addition, positive inter-event times are obtained which guarantee that Zeno behavior does not occur. 
\end{abstract}


\section{Introduction} \label{sec:one}
In recent years, consensus problems over reliable and infinite bandwidth communication networks have been well studied due to their applications in sensor networks and multi-agent systems coordination \cite{Ji07}, \cite{Moreau04}, \cite{Olfati04}. Many consensus protocols rely on the assumption that continuous exchange of information among agents is possible. In many applications, continuous communication is not possible, and it becomes important to discern how frequently agents should communicate to preserve the system properties that stemmed from continuous information exchange. The sampled-data approach has been commonly used to estimate the sampling periods \cite{CaoRen10}, \cite{QinGao12}, \cite{Hayakawa06}, and \cite{Liu10}. An important drawback of periodic transmission is that it requires synchronization between the agents, that is, all agents need to transmit their information at the same time instants and, in some cases, it requires a conservative sampling period for worst case situations. 

More recently, the event-triggered control paradigm has been used to design consensus protocols that account for limited bandwidth communication channels, by reducing the number of transmitted measurements by each agent in the network \cite{Garcia13b}, \cite{kia14}, \cite{Meng13}, \cite{Nowzari14}, \cite{Seyboth13}. Decentralized event-triggered consensus protocols allow agents to decide when to broadcast information independently of each other. 

In event-triggered control we have that measurements are not transmitted periodically in time but they are triggered by the occurrence of certain events. In event-triggered broadcasting  \cite{AntaTabuada10}, \cite{Astrom02}, \cite{Garcia13}, \cite{Donkers12}, \cite{Garcia12IJC}, \cite{Tabuada07}, \cite{WangLemmon11}, \cite{Garcia12CDC}, \cite{Garcia14CDCdis} a subsystem sends its local state to the network only when a measure of the local subsystem state error is above a specified threshold. 
Event-triggered control strategies have also been applied to stabilize multiple coupled subsystems as in \cite{Guinaldo13}, \cite{Guinaldo14}, \cite{Mazo11TAC}, \cite{WangLemmon11}, and \cite{wang2010relaxing}. 

Event-triggered control and communication provides a more efficient method for using network bandwidth. Its implementation in multi-agent systems also provides a highly decentralized way to schedule transmission instants, which eliminates the need for the synchronization required by periodic sampled-data approaches. In this regard, event-triggered control is also robust to clock uncertainties since agents do not need to communicate at the same time instants. 
Different authors have extended the event-triggered consensus approach, for instance,~\cite{ChenHao12} studied event-triggered consensus for discrete time integrators. 
The authors of \cite{YinYue13b} used event-triggered techniques for consensus problems involving a combination of discrete time single and double integrators. 
The authors of \cite{GuoDimarogonas13} studied event-triggered consensus of single integrator systems using nonlinear consensus protocols.
The event-triggered consensus problem with general linear dynamics has been addressed in \cite{Garcia14Auto}, \cite{Garcia14}, \cite{LiuHill12}, \cite{Zhu14} and \cite{Garcia14CDC}. A distributed event-triggered consensus algorithm for a group of agents to determine the global minimizer of a strictly convex cost function was proposed in \cite{kia15}.

In the present paper we consider the event-triggered consensus problem over unreliable communication networks. Previous 
event-triggered consensus approaches such as \cite{Garcia13b}, \cite{Nowzari14}, \cite{Seyboth13} assume that every event-based packet transmission, containing measurement updates, will arrive at its corresponding destination. In many practical scenarios however, this is not the case and packet dropouts occur. We consider the general case of packet dropouts in multi-agent systems where the same transmitted packet may arrive at some destinations and may be lost by other intended receiving agents. We also consider time-varying communication delays. Similarly, delays corresponding to the same broadcasted measurement could be different to every receiving agent.  This type of unreliable communication is \textit{not consistent} with respect to dropouts and delays.

The remainder of this paper is organized as follows. Section~\ref{sec:two} 
describes the problem and the consensus protocol. Section \ref{sec:Results} provides the main results of this paper; asymptotic consensus  
is achieved using a decentralized event-triggered control approach and in the presence of packet losses and delays. Sufficient conditions to guarantee average consensus are investigated in Section~\ref{sec:average}.
Section~\ref{sec:example} presents an illustrative example and Section~\ref{sec:conclusion} concludes the paper.

\section{Preliminaries} \label{sec:two}
\textit{Notation}.
The notations $1_N$ and $0_N$ represent column vectors of all ones and all zeros, respectively. $\mathbb{Z}^{\geq 0}$, $\mathbb{R}$, and $\mathbb{C}$ denote the sets of non-negative integers, real numbers, and complex numbers, respectively. For any $s\in\mathbb{C}$, $Re\left\{s\right\}$ represents the real part of $s$. $J_\nu^{\lambda}$ represents a Jordan block of size $\nu$ corresponding to eigenvalue $\lambda$. The boldface $\textbf{e}^\lambda$ represents the exponential of the scalar $\lambda$ and $\textbf{e}^A$ represents the matrix exponential of matrix $A$. The symbol $\otimes$ denotes the Kronecker product.

\subsection{Graph Theory}
For a team of $N$ agents, the communication among them can be described by a directed graph $\mathcal{G}=\left\{\mathcal{V,E}\right\}$, where $\mathcal{V}=\left\{1,\ldots,N\right\}$ denotes the agent set and $\mathcal{E}\subseteq\mathcal{V}\times\mathcal{V}$ denotes the edge set. An edge $(i,j)$ in the set $\mathcal{E}$ denotes that agent $j$ can obtain information from agent $i$, but not necessarily \emph{vice versa}. For an edge $(i,j)\in \mathcal{E}$, agent $i$ is a neighbor of agent $j$. The set $\mathcal{N}_j$ is called the set of neighbors of agent $j$, and $N_j$ is its cardinality.
A directed path from agent $i$ to agent $j$ is a sequence of edges in a directed graph of the form $(i,p_1),(p_1,p_2),\ldots,(p_{\kappa-1},p_{\kappa})(p_\kappa,j)$, where $p_\ell\in\mathcal{V},~\forall \ell=1,\cdots,\kappa$. A directed graph is \emph{strongly connected} if there is a directed path from every agent to every other agent. A directed graph has \emph{a directed spanning tree} if there exists at least one agent with directed paths to all other agents.

The adjacency matrix $\mathcal{A}\in\mathbb{R}^{N\times N}$ of a directed graph $\mathcal{G}$ is defined by $a_{ij}=1$ if $(j,i)\in\mathcal{E}$ and $a_{ij}=0$ otherwise. The Laplacian matrix $\mathcal{L}$ of $\mathcal{G}$ is defined as $\mathcal{L=D-A}$, where $\mathcal{D}$ represents the degree matrix which is a diagonal matrix with entries $d_{ii}=\sum_{j\in\mathcal{N}_i}a_{ij}$. If a directed graph has a directed spanning tree, then the corresponding Laplacian matrix has only one eigenvalue equal to zero, $\lambda_1=0$, and the following holds for the remaining eigenvalues: $Re\left\{\lambda_i\right\}>0$, for $i=2,...,N$.

\subsection{Problem Statement}
Consider a group of $N$ agents interconnected using a directed communication graph. Each agent can be described by the following: 
\begin{align}
	\dot{x}_i(t)=u_i(t),\ \ i=1,...,N, \label{eq:agents}
\end{align}
 where the local inputs $u_i(t)$ are given by 
\begin{align}
	u_i(t)=-\sum_{j\in\mathcal{N}_i}(x_i(t_{k_i})-x_{ji}(t)),\ \ i=1,...,N, \label{eq:inputs}
\end{align}
where $x_i,u_i\in \mathbb{R}^n$ and $k_i \in \mathbb{Z}^{\geq 0}$ represents the index that defines the non-periodic time sequence of events corresponding to agent $i$, for $i=1,...,N$. Also, $x_{ji}(t)=x_j(t_{k_j^i})$ for $t\in[t_{k_j^i},t_{k_{j}^i+1})$. The variables $x_{ji}(t)$ are piecewise constant and they are not updated at every time $t_{k_j}$ because of packet losses, but only when a measurement from agent $j$ is eventually received by agent $i$ and we denote the time associated with a successful arrival from agent $j$ to agent $i$ as $t_{k_j^i}$. 


We consider the consensus problem over unreliable communication networks while using a decentralized event-triggered control approach. Different challenges are present when transmitting information through a digital communication network. In this paper we jointly consider several issues related to the broadcasting of information in the multi-agent consensus problem.
	 
\begin{enumerate}
	\item The communication channel is of limited bandwidth and, therefore, continuous communication among agents is not possible. A zeno-free, decentralized, event-triggered consensus protocol will be implemented in order to determine the time instants at which each agent should broadcast an information packet containing its current state measurement. 
	\item Time-varying communication delays are present in the communication channel.
 \item Finally, the communication channel may drop information packets so even when events are generated there is no guarantee that all destination agents will receive the transmitted state measurement.
\end{enumerate}

In this paper we consider a general or non-consistent type of delays and packet dropouts. In multi-agent systems consistent delays refer to the case where the delay $d_i(t_{k_i})$ associated with the transmitted state $x_i(t_{k_i})$ is the same for every receiving agent. By non-consistent delays we refer to the more general case where the delay associated to a transmitted state can be different to every receiving agent. In this case we define $d_{ij}(t_{k_i})$ as the time it takes the measurement $x_i(t_{k_i})$ which is released at time $t_{k_i}$ to arrive to agent $j$, for every $j$ such that $i\in\mathcal{N}_j$. For instance, if agents $2$ and $3$ receive information from agent $1$, then the measurement $x_1(t_{k_1})$ released at time $t_{k_1}$ will arrive (if it is not dropped by either agent $2$ or $3$) to agent $2$ at time $t_{k_1}+d_{12}(t_{k_1})$ and to agent $3$ at time $t_{k_1}+d_{13}(t_{k_1})$, where, in general, $d_{12}(t_{k_1})\neq d_{13}(t_{k_1})$. Also, the delay in the communication channel from $i$ to $j$ is time-varying, that is, $d_{ij}(t)$ might not be equal to $d_{ij}(t')$ for $t\neq t'$ and for $i,j=1,...,N$. 

Similarly, we consider non-consistent packet losses \cite{wang2010relaxing}. In \cite{wang2010relaxing} a packet of information  broadcasted by a subsystem is lost in some communication links but it is not lost in other links and it successfully arrives to a subset of nodes. Therefore, some agents may receive different sets of measurements from the same subsystem. 
This means that a broadcasted measurement $x_i(t_{k_i})$ may be successfully received by all, some, or none of the receiving (or destination) agents $j$ such that $i\in\mathcal{N}_j$. 

The main consequence of dealing with non-consistent communication delays and packet dropouts is that agents $j$, for $i\in\mathcal{N}_j$, will hold different versions of the state of agent $i$ since each agent $j$ may receive different updates and also receive them at different time instants. Let $x_{ij}(t)=x_i(t_{k_i^j})$ for $t\in[t_{k_i^j},t_{k_{i}^j+1})$ represent the state $x_i$ as seen by agent $j$, for $i\in\mathcal{N}_j$, where $t_{k_i^j}$ represents the time instants when a measurement from agent $i$ is succesfully received by agent $j$. The state $x_{ij}(t)$ is piecewise constant and it is updated when a measurement transmitted from agent $i$ is successfully received by agent $j$. Define the local state error
\begin{align}
	e_i(t)=x_i(t_{k_i})-x_i(t)   
\end{align}
for $t\in[t_{k_i},t_{k_i+1})$. Also define the state error of agent $i$ as seen by agents $j$, such that $i\in \mathcal{N}_j$,
\begin{align}
	e_{ij}(t)=x_{ij}(t)-x_i(t).  
\end{align}

The error $e_i$ represents the local state error corresponding to agent $i$ and it can be continuously measured in order for agent $i$ to decide when to broadcast its state. On the other hand, the error $e_{ij}$ represents the state error corresponding to agent $i$ as seen by agent $j$ and it cannot be measured by agent $i$ since it does not know the current state $x_{ij}(t)$. Further, the error $e_{ij}$ cannot immediately be reset as the error $e_i$ due to delays and packet losses. 

Assume a uniform Maximum Allowable Number of Successive Dropouts (MANSD) \cite{WangLemmon11}, 
denoted as $\rho-1$, where $\rho>1$ is an integer. This means that if a measurement transmitted by agent $i$ at time $t_{k_i}$ is successfully received by agent $j$, then, at most $\rho-1$ consecutive dropouts are allowed from $i$ to $j$ and, in the worst case, the measurement transmitted at time $t_{k_i+\rho}$ will be successfully received by agent $j$.
In this paper we will determine the largest admissible delay, $d>0$, using the event-triggered communication scheme and based on design parameters. This means that communication delays are time-varying but bounded by $d$, that is, $0\leq d_{ij}(t_{k_i}) \leq d$, for $i,j=1,...,N$. 

\section{Decentralized event-triggered consensus} \label{sec:Results}
In this paper we consider networks of agents that are represented by directed graphs. Also, we assume that the graph has a spanning tree. Let us define agent $r$ as the agent which is the root of the spanning tree. 
 
Under the decentralized event-triggered control protocol, each agent's dynamics can be written as follows
\begin{align}
\left.
	\begin{array}{l l}
	\dot{x}_i(t)&=-\sum_{j\in\mathcal{N}_i,j\neq r}\big(x_i(t)+e_i(t)-x_j(t)-e_{ji}(t)\big)  \\
	  &~~-a_{ir}(x_i(t)+e_i(t)-x_r(t)-e_{ri}(t))  
\end{array}    \right.  \label{eq:agents2}
\end{align}
where $a_{ir}=1$ if agent $r$ is neighbor of agent $i$ and $a_{ir}=0$ otherwise. Define $\eta_i(t)=x_i(t)-x_r(t)$ and we have that
\begin{align}
\left.
	\begin{array}{l l}
	\dot{\eta}_i(t)\!&=\!-\sum_{j\in\mathcal{N}_i,j\neq r}\big(\eta_i(t)-\eta_j(t)\big) \!-\!a_{ir}\eta_i(t)    \\
	  &~~-\sum_{j\in\mathcal{N}_i,j\neq r}\big(e_i(t)-e_{ji}(t)\big) \!-\!a_{ir}(e_i(t)\!-\!e_{ri}(t)) .  
\end{array}    \right.  \label{eq:diffs}
\end{align}
It is convenient to relabel the agents and let $r=1$ to write the Laplacian matrix as follows
\begin{align}
\left.
	\begin{array}{l l}
	\mathcal{L}=\begin{bmatrix}
	 0 &0_{N-1}^\text{T} \\
	-a_1 & \mathcal{L}_r
\end{bmatrix}.
\end{array}    \right.  
\end{align}
where $a_1=[a_{21} \ a_{31} \ldots a_{N1}]^\text{T}$. Define $\eta=[\eta_2^\text{T} \ \eta_3^\text{T}\ldots \eta_N^\text{T}]^\text{T}$ and 
we can write the overall system as follows
\begin{align}
	\dot{\eta}(t)=-(\mathcal{L}_r\otimes I_n)\eta(t)+ \xi(t)   \label{eq:OA}
\end{align}
where 
\begin{align}
	\xi(t)=
\begin{bmatrix}
	 \sum_{j\in\mathcal{N}_2}(e_{j2}(t)-e_2(t)) \\
	\sum_{j\in\mathcal{N}_3}(e_{j2}(t)-e_3(t))  \\
	\vdots  \\
	\sum_{j\in\mathcal{N}_N}(e_{jN}(t)-e_N(t))
\end{bmatrix}.   \label{eq:xi}
\end{align}
Note that the eigenvalues of $\mathcal{L}_r$ are $\lambda_i$, for $i=2,...,N$, that is, the eigenvalues of $\mathcal{L}$ with positive real parts. Therefore, there exist $\hat{\beta}, \hat{\lambda}>0$ such that the relation $ \left\|\textbf{e}^{-\mathcal{L}_rt}\right\|\leq\hat{\beta}\textbf{e}^{-\hat{\lambda} t}$ holds. 

In the particular case that the graph is strongly connected any agent may serve as agent $r=1$. Eq. \eqref{eq:diffs} is now given by
\begin{align}
\left.
	\begin{array}{l l}
	\dot{\eta}_i(t)\!&=\!-\sum_{j\in\mathcal{N}_i,j\neq r}\big(\eta_i(t)-\eta_j(t)\big) \!-\!a_{ir}\eta_i(t)    \\
	  &~~-\sum_{j\in\mathcal{N}_i,j\neq r}\big(e_i(t)-e_{ji}(t)\big) \!-\!a_{ir}(e_i(t)\!-\!e_{ri}(t))  \\
		&~~+\sum_{j\in\mathcal{N}_r}\big(-\eta_j(t)+e_r(t)-e_{jr}(t)\big)   
\end{array}    \right.  
\end{align}
and the overall system is given by 
\begin{align}
	\dot{\eta}(t)=-(\mathcal{L}_s\otimes I_n)\eta(t)+ \xi(t) -1_{N-1}\otimes \xi_r(t)  
\end{align}
where $\mathcal{L}_s=\mathcal{L}_r+1_{N-1}\alpha^\text{T}$, $\alpha=[a_{12} \ a_{13}\ldots a_{1N}]^\text{T}$, and $\xi_r(t)=\sum_{j\in\mathcal{N}_r}(e_{jr}(t)-e_r(t))$. 
Let 
\begin{align}
	S=\begin{bmatrix}
    1 & 0^\text{T}_{N-1} \\
	 1_{N-1} & I_{N-1}
\end{bmatrix}  \nonumber
\end{align}
and compute 
\begin{align}
	\hat{\mathcal{L}}=S^{-1}\mathcal{L}S=\begin{bmatrix}
    0 & -\alpha^\text{T} \\
	 0_{N-1} & \mathcal{L}_s
\end{bmatrix} . \nonumber
\end{align}
Hence, the eigenvlaues of $\mathcal{L}_s$ are $\lambda_i$, for $i=2,...,N$.

Due to presence of packet losses each agent will impose a maximum time between events $\bar{\delta}>0$, that is, an event is generated if $t=t_{k_i}+ \bar{\delta}$. The need for this maximum inter-event time is explained at the end of this section. 
The following theorem establishes asymptotic convergence of the event-triggered consensus algorithm in the presence of packet dropouts and delays. It also provides a positive lower-bound on the inter-event time intervals and an estimate of the largest admissible delay for given value of the MANSD.

\begin{theorem}\label{th:consensus}
Assume that the communication graph has a spanning tree and the MANSD is $\rho-1$, for $\rho>1$. Then, for some $\gamma_d>0$, agents~\eqref{eq:agents} implementing decentralized control inputs~\eqref{eq:inputs} achieve consensus asymptotically in the presence of communication delays $d_{ij}(t_{k_i})\leq d$, 
if agent $i$'s events, $t_{k_i}$ for i=1,\dots,N, are generated according to the following condition
\begin{align}
	   t_{k_i+1} \!\!=\! \arg\min \!\left\{t>t_{k_i} \big| \left\|e_i(t)\right\|= \beta \normalfont{\textbf{e}}^{-\lambda t} \text{or} \ t\!=\!t_{k_i}\!\! +\!\bar{\delta} \right\}  \label{eq:thre}
\end{align}
where $t_{0}=0$, $k_i\in \mathbb{Z}^{\geq 0}$, $\beta>0$ and $0<\lambda<\hat{\lambda}$, and 
\begin{align}
    \left.
	\begin{array}{l l}
	d=\frac{1}{\hat{\lambda}}\ln\Big(1 + \frac{\gamma_d}{\frac{H_2}{\lambda} + \frac{H_1}{\hat{\lambda}}\normalfont{\textbf{e}}^{(\lambda-\hat{\lambda}) t_{k_i+\rho}}}\Big).
	\end{array}  \label{eq:delay}   \right.
\end{align}
Furthermore, the agents do not exhibit Zeno behavior and the inter-event times $t_{k_i+1}-t_{k_i}$ for every agent $i=1,...,N$ are bounded below by the \textit{positive} time $\tau$, that is
\begin{align}
  \tau\leq t_{k_i+1}-t_{k_i}  \label{eq:tk}
\end{align}	
where 
\begin{align}
    \left.
	\begin{array}{l l}
	\tau=\frac{1}{\hat{\lambda}}\ln\Big(1 + \frac{\beta}{\frac{H_2}{\lambda} + \frac{H_1}{\hat{\lambda}}\normalfont{\textbf{e}}^{(\lambda-\hat{\lambda}) t_{k_i}}}\Big)
	\end{array}  \label{eq:tau}   \right.
\end{align}
and $H_1$ and $H_2$ are given by \eqref{eq:Hs} below.
\end{theorem}

\textit{Proof}.
To prove Theorem \ref{th:consensus}, the following observations are required.  Note that because of threshold \eqref{eq:thre}, the error $e_i$ is reset to zero at the event instants $t_{k_i}$, that is, $e_i(t_{k_i})=0$. Thus, the error $e_i$ satisfies $\left\|e_i(t)\right\|\leq \beta \textbf{e}^{-\lambda t}$ and we have that $\left\|e(t)\right\|\leq\sqrt{N} \beta \textbf{e}^{-\lambda t}$, where $e(t)=[e_1(t)^\text{T} \ e_2(t)^\text{T}\ldots e_N(t)^\text{T}]^\text{T}$.

\textit{Analyze packet losses:}
Let us consider first the part of the state errors $e_{ij}$ due to packet losses. Assume without loss of generality that the last  update transmitted by agent $i$ and successfully received by agent $j$ takes place at time $t_{k_i}$. This means that $t_{k_i^j}=t_{k_i}$. Agent $i$ will generate the next event at time $t_{k_i+1}$ and it will broadcast its current state $x_i(t_{k_i+1})$. This event is generated because the condition \eqref{eq:thre} is satisfied, then, we have that $\left\|e_i(t_{k_i+1})\right\|=\left\|x_i(t_{k_i})-x_i(t_{k_i+1})\right\| \leq \beta\textbf{e}^{-\lambda t_{k_i+1}}$. Assume that the update at time $t_{k_i+1}$ is dropped, so $\left\|e_{ij}(t_{k_i+1})\right\| \leq \beta\textbf{e}^{-\lambda t_{k_i+1}}$. Similarly, agent $i$ will generate the following event at time $t_{k_i+2}$ (and broadcast its current state $x(t_{k_i+2})$) and we have that $\left\|e_i(t_{k_i+2})\right\| \leq \beta\textbf{e}^{-\lambda t_{k_i+2}}$ is satisfied. If the number of successive dropouts after the last successful update is the MANSD, $\rho -1$, then we have that the local error at time $t_{k_i+\rho}$, just before the update $x(t_{k_i+\rho})$ is successfully received by agent $j$, is
\begin{align}
  \left.
	\begin{array}{l l}
  \left\|e_i(t_{k_i+\rho})\right\|\!\!&=\left\|x_i(t_{k_i+\rho-1})-x_i(t_{k_i+\rho})\right\|  \\
	  &=\left\|x_i(t_{k_i+\rho-1})-x_i(t_{k_i+\rho}) \right. \\
	  &~~+x_i(t_{k_i})-x_i(t_{k_i}) \ldots  \\ 
		&~~\left.+x_i(t_{k_i+\rho-2})-x_i(t_{k_i+\rho-2}) \right\|  \\
		&=\left\|	x_i(t_{k_i})-x_i(t_{k_i+\rho}) - e_i(t_{k_i+1})  
		- e_i(t_{k_i+2}) - \ldots -e_i(t_{k_i+\rho-1})\right\| \\
		&\geq 	\left\|x_i(t_{k_i})-x_i(t_{k_i+\rho})\right\| -\beta\textbf{e}^{-\lambda t_{k_i+1}} 
		 - \beta\textbf{e}^{-\lambda t_{k_i+2}} -\ldots -\beta\textbf{e}^{-\lambda t_{k_i+\rho-1}}.
\end{array}   \label{eq:drops2} \right.
\end{align}
Thus, the error $e_{ij}$ when the MANSD occurs is bounded by 
\begin{align}
  \left.
	\begin{array}{l l}
    \left\|e_{ij}(t_{k_i+\rho})\right\|&= \left\|x_i(t_{k_i})-x_i(t_{k_i+\rho})\right\|  \\
		 & \leq \beta\big(\textbf{e}^{-\lambda t_{k_i+1}} + \ldots + \textbf{e}^{-\lambda t_{k_i+\rho}}   \big)
\end{array}   \label{eq:drops3} \right.
\end{align}
Let $t=t_{k_i}+s$, for $s\in[0,s_\rho)$. Thus, the time instant $t=t_{k_i+1}$ can be represented by $t=t_{k_i}+s_1$. Similarly, $t_{k_i+2}=t_{k_i}+s_2$ and so on. Note that the event at time $t_{k_i+\mu}$ satisfies the following: $\mu\tau\leq t_{k_i+\mu}-t_{k_i}\leq \mu\bar{\delta}$   , for $\mu=1,...,\rho$, where $\tau\geq0$ represents the minimum inter-event time, that is, $\mu\tau\leq s_\mu\leq \mu\bar{\delta}$. We can write
\begin{align}
  \left.
	\begin{array}{l l}
   \left\|e_{ij}(t_{k_i+\rho})\right\| &\leq \beta\textbf{e}^{-\lambda t_{k_i}}\big(\textbf{e}^{-\lambda s_1} + \ldots + \textbf{e}^{-\lambda s_\rho} \big)  \\
	&\leq \beta\textbf{e}^{-\lambda t_{k_i}} \sum_{\mu=1}^\rho \textbf{e}^{-\mu\lambda\tau}.
\end{array}   \label{eq:drops4} \right.
\end{align}
We also have that $\left\|e_{ij}(t)\right\| \leq \beta\textbf{e}^{-\lambda t_{k_i}} \sum_{\mu=1}^\rho \textbf{e}^{-\mu\lambda\tau}$, for $t\in[t_{k_i},t_{k_i+\rho})$, since $\rho$ represents the worst case number of consecutive packet losses. This can be simply shown by letting $\rho'<\rho$ and if some packet is received at $t_{k_i+\rho'}$ then we would have $\sum_{\mu=1}^{\rho'} \textbf{e}^{-\mu\lambda\tau}\leq\sum_{\mu=1}^\rho \textbf{e}^{-\mu\lambda\tau}$ since $\rho'<\rho$. Then, in general, we have that
\begin{align}
  \left.
	\begin{array}{l l}
   \left\|e_{ij}(t)\right\|	&\leq \beta\textbf{e}^{\lambda s} \sum_{\mu=1}^\rho \textbf{e}^{-\mu\lambda\tau}\textbf{e}^{-\lambda t}  
	 	\leq \gamma_l\textbf{e}^{-\lambda t}  
\end{array}   \label{eq:drops5} \right.
\end{align}
for $t\in[t_{k_i},t_{k_i+\rho})$, where
\begin{align}
\gamma_l=\beta\textbf{e}^{\rho\lambda \bar{\delta}} \sum_{\mu=1}^\rho \textbf{e}^{-\mu\lambda\tau}.
\end{align}
Note that in the derivation of $\gamma_l$ the relation $\tau\geq 0$ was used. Although this relation is sufficient to obtain $\gamma_l$, we will show later in this proof that $\tau>0$ in order to avoid Zeno behavior. 

\begin{figure}
	\begin{center}
		\includegraphics[width=12cm,height=6.5cm,trim=1.4cm 3.8cm 2.4cm .6cm]{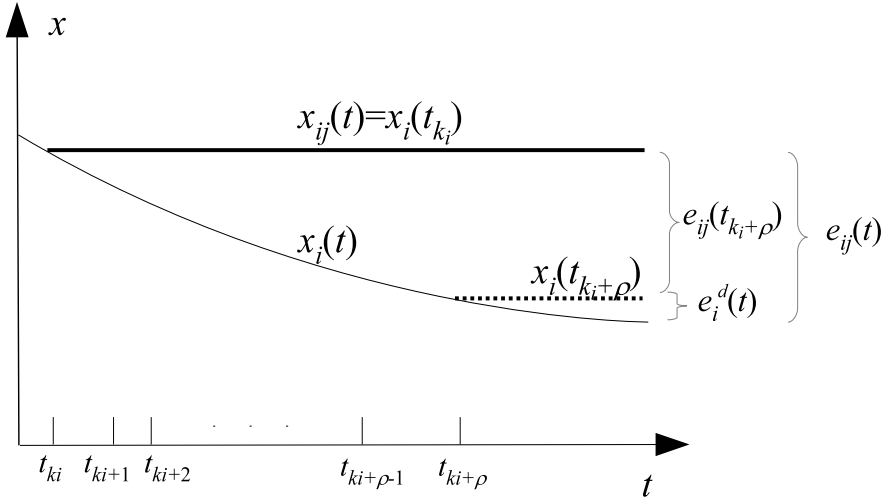}
	\caption{State errors due to packet dropouts and delays during the time interval $t\in[t_{k_i+\rho},t_{k_i+\rho}+d)$}
	\label{fig:Error}
	\end{center}
\end{figure}

Due to communication delays, the update $x_i(t_{k_i+\rho})$ will not be received by agent $j$ until time $t_{k_i+\rho}+d$, in the worst case delay. Hence, we bound the error $e_{ij}(t)$ within the extended time interval $t\in[t_{k_i},t_{k_i+\rho}+d)$.  
Define the error due to delays as follows 
\begin{align}
   e_{i}^d(t)=x_i(t_{k_i+\rho})-x_i(t)
\end{align}
for $t\in[t_{k_i+\rho},t_{k_i+\rho}+d)$. Let us write the following
\begin{align}
    \left.
	\begin{array}{l l}
	  e_{ij}(t)&=x_i(t_{k_i})-x_i(t)  \\
	      &= x_i(t_{k_i})- x_i(t_{k_i+\rho})+e_i^d(t).
	\end{array}    \right.
\end{align}
Thus, the following holds
\begin{align}
  \left\|e_{ij}(t)\right\|&\leq  \gamma_l\textbf{e}^{-\lambda t}  + \left\|e_i^d(t)\right\|  \label{eq:lderr}
\end{align}
These relationships are illustrated in Fig. \ref{fig:Error} for a single element of the involved variables. Given a $\gamma_d>0$ we can always guarantee that there exist some $d\geq 0$ such that $\left\|e_i^d(t)\right\|\leq \gamma_d \textbf{e}^{-\lambda t}$ for $t\in[t_{k_i+\rho},t_{k_i+\rho}+d)$. The previous relation holds because $\gamma_d>0$, $e_i^d(t_{k_i+\rho})=0$, and  $e_i^d$ is continuous in the interval $t\in[t_{k_i+\rho},t_{k_i+\rho}+d)$.
However, in the last step of this proof, an estimate of the largest admissible delay will be obtained as a function of the design parameter $\gamma_d$. At this point we have obtained the following relationship involving $e_{ij}(t)$ 
\begin{align}
  \left\|e_{ij}(t)\right\| \leq \gamma \textbf{e}^{-\lambda t}  \label{eq:errijb}
\end{align}
for $i,j=1,...,N$, where $\gamma=\gamma_l+\gamma_d$.

In the remainder of the proof, it is shown that agents~\eqref{eq:agents} with decentralized control inputs~\eqref{eq:inputs} converge asymptotically towards consensus. Also, the lower bound, $\tau$, on inter-event times for every agent is obtained. Finally, given the parameter $\gamma_d,$ the worst case admissible delay is determined.

\textit{Convergence:} Let us write \eqref{eq:xi} as  $\xi=\xi_d-(\mathcal{D}\otimes I_n)e$, where 
\begin{align}
	\xi_d(t)=
\begin{bmatrix}
	 \sum_{j\in\mathcal{N}_1}e_{j1}(t) \\
	\sum_{j\in\mathcal{N}_2}e_{j2}(t)  \\
	\vdots  \\
	\sum_{j\in\mathcal{N}_N}e_{jN}(t)
\end{bmatrix}.   \label{eq:xid}
\end{align}
Then we have that
\begin{align}
  \left\|\xi(t)\right\| \leq \Gamma \textbf{e}^{-\lambda t}  \label{eq:errxi}
\end{align}
where $\Gamma=\bar{N}_i\sqrt{N}\beta + \gamma \big(\sum_{i=1}^N N_i^2\big)^{1/2}$ and $\bar{N}_i=\max_iN_i$.

From \eqref{eq:OA}, the response of $\eta(t)$ satisfies
\begin{align}
  \left.
	\begin{array}{l l}
  \left\|\eta(t)\right\|&=\left\|\textbf{e}^{-\mathcal{L}_rt}\eta(0) + \int_0^t \textbf{e}^{-\mathcal{L}_r(t-s)}\xi(s)ds \right\|   \\
	   &\leq \hat{\beta}\eta_0\textbf{e}^{-\hat{\lambda} t} + \hat{\beta}\Gamma \int_0^t \textbf{e}^{-\hat{\lambda}(t-s)}\textbf{e}^{-\lambda s} ds  \\
    &\leq \hat{\beta}\eta_0\textbf{e}^{-\hat{\lambda} t} + \frac{\hat{\beta}\Gamma}{\hat{\lambda}-\lambda} (\textbf{e}^{-\lambda t}-\textbf{e}^{-\hat{\lambda} t}) 
\end{array}   \nonumber \right.
\end{align}
where $\eta_0=\left\|\eta(0)\right\|$. As time goes to infinity we have that
\begin{align}
\lim_{t\rightarrow\infty}  \left\|\eta(t)\right\| =0.
\end{align}
and the agents achieve consensus asymptotically. 

\begin{figure}
	\begin{center}
		\includegraphics[width=12cm,height=7.5cm,trim=.9cm .7cm .9cm .3cm]{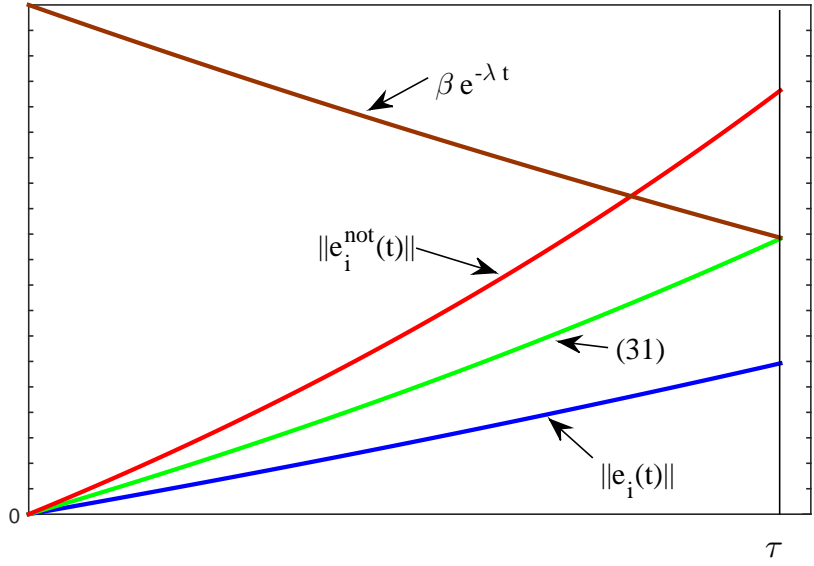}
	\caption{Growth of error $\left\|e_i(t)\right\|$ is bounded by \eqref{eq:errresp}}
	\label{fig:exps}
	\end{center}
\end{figure}

\textit{Determine $\tau$:} In order to establish a positive lower bound on the inter-event time intervals (and avoid Zeno behavior) we look at the dynamics of the error $e_i(t)$, $\frac{d}{dt}e_i(t)=\dot{e}_i(t)=-\dot{x}_i(t)$ for $t\in[t_{k_i},t_{k_i+1})$ with $e_i(t_{k_i})=0$. We have
\begin{align}
    \left.
	\begin{array}{l l}
	  \left\|\dot{e}_i(t)\right\|& \leq \left\|\sum_{j\in\mathcal{N}_i,j\neq r}\big(\eta_i(t)-\eta_j(t)\big) +a_{ir}\eta_i(t) \right\|   
	  +\left\|\sum_{j\in\mathcal{N}_i}\big(e_i(t)-e_{ji}(t)\big) \right\|       \\
			&\leq \left\|(\mathcal{L}_r\otimes I_n)\eta(t)\right\| \!+\! N_i\beta\textbf{e}^{-\lambda t}\!+\!N_i\gamma\textbf{e}^{-\lambda t}.
	\end{array}     \right.  \label{eq:errdyn}
\end{align}
Let $L=\left\|\mathcal{L}_r\right\|$, then
\begin{align}
    \left.
	\begin{array}{l l}
	  \left\|\dot{e}_i(t)\right\|&\leq L\Big( \hat{\beta}\eta_0\textbf{e}^{-\hat{\lambda} t} + \frac{\hat{\beta}\Gamma}{\hat{\lambda}-\lambda} (\textbf{e}^{-\lambda t}-\textbf{e}^{-\hat{\lambda} t}) \Big) 
		  + N_i(\beta+\gamma)\textbf{e}^{-\lambda t}.
	\end{array}  \label{eq:errdyn2}   \right.
\end{align}
Define
\begin{align}
    \left.
	\begin{array}{l l}
  H_1=L\hat{\beta}\Big(\eta_0 - \frac{\Gamma}{\hat{\lambda}-\lambda}\Big)  \\
  H_2=N_i(\beta+\gamma) + \frac{L\hat{\beta}\Gamma}{\hat{\lambda}-\lambda}.
	\end{array}  \label{eq:Hs}   \right.
\end{align}
The error response during the time interval $t\in[t_{k_i},t_{k_i+1})$ can be bounded as follows
\begin{align}
    \left.
	\begin{array}{l l}
	  \left\|e_i(t)\right\|&\leq \int_{t_{k_i}}^t\big(H_1\textbf{e}^{-\hat{\lambda} s} + H_2\textbf{e}^{-\lambda s} \big)ds \\
	&\leq \frac{H_1}{\hat{\lambda}}\big(\textbf{e}^{-\hat{\lambda} t_{k_i}}-\textbf{e}^{-\hat{\lambda} t}\big)  
	+ \frac{H_2}{\lambda}\big(\textbf{e}^{-\lambda t_{k_i}}-\textbf{e}^{-\lambda t}\big)  \\
	&\leq \frac{H_1}{\hat{\lambda}}\big(1-\textbf{e}^{-\hat{\lambda} \tau}\big)\textbf{e}^{-\hat{\lambda} t_{k_i}} 
	+ \frac{H_2}{\lambda}\big(1-\textbf{e}^{-\lambda \tau}\big)\textbf{e}^{-\lambda t_{k_i}} 
	\end{array}  \label{eq:errresp}   \right.
\end{align}
where $\tau=t-t_{k_i}$. 
Thus, the time $\tau>0$ that it takes for the last expression in \eqref{eq:errresp} to grow from zero, at time $t_{k_i}$, to reach the threshold $\beta\textbf{e}^{-\lambda t}=\beta\textbf{e}^{-\lambda(t_{k_i}+\tau)}$ is less or equal than the time it takes the error $\left\|e_i(t)\right\|$ to grow from zero, at time $t_{k_i}$, to reach the same threshold and generate the following event at time $t_{k_i+1}$. This means that $0<\tau\leq t_{k_i+1}-t_{k_i}$. 
These relationships are shown in Fig. \ref{fig:exps}. We know that $\left\|e_i(t)\right\|$ will be bounded by \eqref{eq:errresp}, for instance, it will grow as shown by the curve in blue in Fig. \ref{fig:exps}. In the worst case, it will grow as the right hand term in \eqref{eq:errresp} (shown in green in the same figure), when \eqref{eq:errresp} is satisfied with equality. However, we will never have the case $\left\|e_i^{not}(t)\right\|$  shown by the curve in red. Hence, \eqref{eq:errresp} will hit the threshold  $\beta \textbf{e}^{-\lambda t}=\beta \textbf{e}^{-\lambda (t_{k_i}+\tau)}$ sooner (at the same time in the worst case) than $\left\|e_i(t)\right\|$.
Thus, we wish to find a lower-bound $\tau>0$ such that the following holds
\begin{align}
    \left.
	\begin{array}{l l}
	\frac{H_1}{\hat{\lambda}}\big(1-\textbf{e}^{-\hat{\lambda} \tau}\big)\textbf{e}^{-\hat{\lambda} t_{k_i}} 
	+ \frac{H_2}{\lambda}\big(1-\textbf{e}^{-\lambda \tau}\big)\textbf{e}^{-\lambda t_{k_i}} \leq \beta\textbf{e}^{-\lambda (t_{k_i}+\tau)}
	\end{array}  \label{eq:errresp2}   \right.
\end{align}
equivalently,
\begin{align}
    \left.
	\begin{array}{l l}
	\frac{H_1}{\hat{\lambda}}\big(1\!-\!\textbf{e}^{-\hat{\lambda} \tau}\big)\textbf{e}^{(\lambda-\hat{\lambda}) t_{k_i}} 
	\!+\! \frac{H_2}{\lambda}\big(1\!-\!\textbf{e}^{-\lambda \tau}\big) \leq \beta\textbf{e}^{-\lambda \tau}.
	\end{array}  \label{eq:errresp3}   \right.
\end{align}
Inequality \eqref{eq:errresp3} can also be written as follows
\begin{align}
    \left.
	\begin{array}{l l}
	\frac{H_1}{\hat{\lambda}}\textbf{e}^{(\lambda-\hat{\lambda}) t_{k_i}} + \frac{H_2}{\lambda}  
	\leq \beta\textbf{e}^{-\lambda \tau} + \frac{H_2}{\lambda}\textbf{e}^{-\lambda \tau} + \frac{H_1}{\hat{\lambda}}\textbf{e}^{-\hat{\lambda} \tau}\textbf{e}^{(\lambda-\hat{\lambda}) t_{k_i}} 
	\end{array}  \label{eq:errresp4}   \right.
\end{align}
because $\textbf{e}^{-\lambda \tau}\geq \textbf{e}^{-\hat{\lambda} \tau}$ for $\hat{\lambda}>\lambda$ and for any $\tau\geq 0$. We also have that
\begin{align}
    \left.
	\begin{array}{l l}
	 \beta\textbf{e}^{-\lambda \tau} + \frac{H_2}{\lambda}\textbf{e}^{-\lambda \tau} + \frac{H_1}{\hat{\lambda}}\textbf{e}^{-\hat{\lambda} \tau}\textbf{e}^{(\lambda-\hat{\lambda}) t_{k_i}} 
	\geq \big(\beta + \frac{H_2}{\lambda} + \frac{H_1}{\hat{\lambda}}\textbf{e}^{(\lambda-\hat{\lambda}) t_{k_i}}\big)   \textbf{e}^{-\hat{\lambda} \tau}
	\end{array}  \label{eq:errresp4-1}   \right.
\end{align}
Here, we strive to obtain the largest value of $\tau$ such that \eqref{eq:errresp4} holds. This value can be obtained by solving the equation 
\begin{align}
    \left.
	\begin{array}{l l}
	\frac{H_1}{\hat{\lambda}}\textbf{e}^{(\lambda-\hat{\lambda}) t_{k_i}} + \frac{H_2}{\lambda}  
	=\big(\beta + \frac{H_2}{\lambda} + \frac{H_1}{\hat{\lambda}}\textbf{e}^{(\lambda-\hat{\lambda}) t_{k_i}}\big)   \textbf{e}^{-\hat{\lambda} \tau}
	\end{array}  \nonumber  \right.
\end{align}
we guarantee that \eqref{eq:errresp2} holds.
The explicit solution for the lower bound on the inter-event time intervals is given by \eqref{eq:tau}.
By the selection $\hat{\lambda}>\lambda$, we have that $\textbf{e}^{(\lambda-\hat{\lambda})t_{k_i}}\leq 1$ for any $t_{k_i}\geq 0$, ensuring that $\tau>0$.

\textit{Worst case delay:} Let us now estimate the largest admissible communication delay $d$ subject to the constraint $\left\|e_i^d(t)\right\|\leq \gamma_d\textbf{e}^{-\lambda t}$. 
We have that the part of the state error due to communication delays satisfies $\dot{e}_i^d=-\dot{x}_i$ for $t\in[t_{k_i+\rho},t_{k_i+\rho}+d)$, with $e_i^d(t_{k_i+\rho})=0$. Then, we can follow similar steps as in \eqref{eq:errdyn}-\eqref{eq:errresp4} to show first that the error $e_i^d$ satisfies the following
\begin{align}
    \left.
	\begin{array}{l l}
	  \left\|e_i^d(t)\right\|\!	&\leq\! \frac{H_1}{\hat{\lambda}}\big(1\!-\!\textbf{e}^{-\hat{\lambda} d}\big)\textbf{e}^{-\hat{\lambda} t_{k_i+\rho}} 
		\!+\! \frac{H_2}{\lambda}\big(1\!-\!\textbf{e}^{-\lambda d}\big)\textbf{e}^{-\lambda t_{k_i+\rho}} 
	\end{array}   \right. \nonumber 
\end{align}
for $t\in[t_{k_i+\rho},t_{k_i+\rho}+d)$. Thus, in order to guarantee that $\left\|e_i^d(t)\right\|\leq \gamma_d \textbf{e}^{-\lambda t}$ we solve for $d$ in the following
\begin{align}
    \left.
	\begin{array}{l l}
	\frac{H_1}{\hat{\lambda}}\big(1-\textbf{e}^{-\hat{\lambda} d}\big)\textbf{e}^{-\hat{\lambda} t_{k_i+\rho}} 
	+ \frac{H_2}{\lambda}\big(1-\textbf{e}^{-\lambda d}\big)\textbf{e}^{-\lambda t_{k_i+\rho}} \leq \gamma_d\textbf{e}^{-\lambda (t_{k_i+\rho}+d)}.
	\end{array}  \label{eq:delerrp2}   \right.
\end{align}
An explicit solution for $d$ is given by \eqref{eq:delay}. For any communication delay $d_{ij}\leq d$ the inequality \eqref{eq:delerrp2} holds and the proof is complete. $\square$

\textit{Remark 1}. The maximum inter-event time which was included as the second condition in \eqref{eq:thre}, guarantees that the sequence $t_{k_i}\rightarrow\infty$. In the presence of packet dropouts and using a protocol only based on error events (first condition in \eqref{eq:thre}) the following situation may appear. The local state error of some agent may satisfy $e_i(t)=0$, for $t\geq t_{k_i}$. This means that after some update at time $t_{k_i}$, agent $i$'s real state does not change and no further events are triggered; however, since packet dropouts are possible, some or all agents that receive information from agent $i$ may never receive the final update. The maximum inter-event time ensures that all agents are finally updated with the correct state value even if the error remains equal to zero after some update. Note that the time event is not needed if packet dropouts do not exist and every measurement update always arrives at its destination.

\textit{Remark 2}. As with many consensus algorithms, an estimate of the second smallest eigenvalue of the Laplacian matrix is required in order to estimate $\hat{\lambda}$; this is the only global information needed by the agents.  Algorithms for distributed estimation of the second eigenvalue of the Laplacian have been presented in \cite{Aragues12}, and \cite{Franceschelli09}. Readers are referred to these papers for details.

\textit{Remark 3}. Note that agents do not need to use the same threshold parameters $\beta$ and $\lambda$. Same results are obtained when every agent implements the event rule in \eqref{eq:thre} using its own $\beta_i$ and $\lambda_i$. All results in Theorem \ref{th:consensus} can be obtained by defining $\beta=\max_i \beta_i$ and $\lambda=\min_i \lambda_i$ and noting that $\left\|e_i(t)\right\|\leq \beta_i\textbf{e}^{-\lambda_i t}\leq \beta\textbf{e}^{-\lambda t}$. Thus, each agent only needs an estimate $\lambda_i$ of $\hat{\lambda}$ such that $\lambda_i\leq \hat{\lambda}$ and the agents do not need to synchronize what threshold parameters to implement.

\textit{Remark 4}. Another practical and decentralized characteristic of the event-triggered consensus protocol discussed in this paper is that it allows for synchronization errors in the agents clocks. In order to implement this event-triggered strategy we need to evaluate the function $\textbf{e}^{-\lambda t}$. However, it is not necessary for the time index $t$ to be the same for every agent. Consider the local threshold function for agent $i$ to be $\textbf{e}^{-\lambda t^i}$, where, $t^i$ represents the time as given by its local clock. Then, there exist a $t$ such that $t^i=t+t_\delta^i$, where $t_\delta^i\geq 0$ represents the clock offset of agent $i$ with respect to $t$. Given the previous relation and the fact that the local events for agent $i$ are triggered using the function $\textbf{e}^{-\lambda t^i}$, we have that
\begin{align}
   \left\|e_i(t)\right\| \leq \beta\textbf{e}^{-\lambda t^i} \leq \beta\textbf{e}^{-\lambda t}.
\end{align}
Similarly, we also obtain 
\begin{align}
   \left\|e_{ij}(t)\right\| \leq \gamma\textbf{e}^{-\lambda t^i} \leq \gamma\textbf{e}^{-\lambda t}.
\end{align}
Based on these two relationships, we can show that the same properties of the consensus protocol hold by following the remaining steps of the proof of Theorem \ref{th:consensus}.

\section{Average consensus} \label{sec:average} 
It was shown in Section \ref{sec:Results} that agents can reach agreement asymptotically using a decentralized event-triggered consensus protocol in the presence of non-consistent packet dropouts and non-consistent communication delays. However, under these conditions, agents do not reach average consensus in general, not even if the communication graph is undirected. 

In the case where delays and packet losses are consistent (same for every receiving agent), a modified consensus protocol can be implemented in order to obtain average consensus. The main difference is the implementation of acknowledgment (ACK) messages that are used when a packet is successfully received by the destination agents. The local control inputs are now given by
\begin{align}
	u_i(t)=-\sum_{j\in\mathcal{N}_i}(x_{ij}(t)-x_{ji}(t)),\ \ i=1,...,N. \label{eq:inputsco}
\end{align}

The following corollary establishes asymptotic convergence of the event-triggered consensus algorithm to the initial average in the presence of consistent delays and packet dropouts. 
\begin{corollary}  \label{co:avg}
Assume that the communication graph is undirected and connected and the MANSD is $\rho>1$. Also assume that ACK messages are transmitted without delay and are never dropped. Then, for some $\gamma_d>0$, agents~\eqref{eq:agents} with decentralized control inputs~\eqref{eq:inputsco} achieve average consensus asymptotically in the presence of constant communication delays $d_i\leq d$ 
if the events are generated according to \eqref{eq:thre}, for $i=1,...,N$, where $\beta>0$ and $0<\lambda<\hat{\lambda}$, and 
\begin{align}
    \left.
	\begin{array}{l l}
	d=\frac{1}{\hat{\lambda}}\ln\Big(1+ \frac{\gamma_d}{\frac{K_2}{\lambda} + \frac{K_1}{\hat{\lambda}}\normalfont{\textbf{e}}^{(\lambda-\hat{\lambda}) t_{k_i+\rho}}}\Big)
	\end{array}  \label{eq:delayco}   \right.
\end{align}
Furthermore, the agents do not exhibit Zeno behavior and the inter-event times $t_{k_i+1}-t_{k_i}$ for every agent $i=1,...,N$ are bounded by the \textit{positive} time $\tau$ as in \eqref{eq:tk} where 
\begin{align}
    \left.
	\begin{array}{l l}
	\tau=\frac{1}{\hat{\lambda}}\ln\Big(1+ \frac{\beta}{\frac{K_2}{\lambda} + \frac{K_1}{\hat{\lambda}}\normalfont{\textbf{e}}^{(\lambda-\hat{\lambda}) t_{k_i}}}\Big)
	\end{array}  \label{eq:tauco}   \right.
\end{align}
and
\begin{align}
    \left.
	\begin{array}{l l}
  K_1=L\hat{\beta}\Big(\eta_0 - \frac{\gamma\sqrt{N}L'}{\hat{\lambda}-\lambda}\Big)  \\
  K_2=2N_i\gamma + \frac{L\hat{\beta}\gamma\sqrt{N}L'}{\hat{\lambda}-\lambda}.
	\end{array}  \label{eq:Ks}   \right.
\end{align}
$L'=\left\|\mathcal{L}'\right\|$, $\mathcal{L}'=\mathcal{L}_{2:N}-1_{N-1}\mathcal{L}_1$, and  $\mathcal{L}_{2:N}$ represents the matrix obtained by removing the first row, $\mathcal{L}_1$, of $\mathcal{L}$.
\end{corollary}
\textit{Proof}. 
 In this case agent $i$ keeps two sampled versions of its own state: $x_i(t_{k_i})$ and $x_{ij}(t)$. The first variable is used to compute new events and is updated at every local time event $t_{k_i}$. The second sample is used to compute its own control input as in \eqref{eq:inputsco} and it is updated only when an ACK message is received from agents $j$, $i \in \mathcal{N}_j$. The end result of this protocol is that every agent $j$, $i \in \mathcal{N}_j$, has the same version $x_{ij}(t)$ of the state of agent $i$.

We define $e_{d_i}(t)=x_{ij}(t)-x_i(t)$ define $e_d=[e_{d_1}^\text{T} \ e_{d_2}^\text{T} \cdots e_{d_N}^\text{T}]^\text{T}$. Here, we do not use $e_{ij}$ since $x_{ij}$ is the same for every agent $j \in \mathcal{N}_i$.
The same steps in \eqref{eq:drops2}-\eqref{eq:errijb} can be followed to obtain 
\begin{align}
  \left\|e_{d_i}(t)\right\| \leq \gamma \textbf{e}^{-\lambda t}  \nonumber
\end{align}
since the corresponding packet dropout analysis in Theorem \ref{th:consensus} holds for any pair of transmitting agent $i$ and receiving agent $j$.

We can write $\dot{x}(t)=-(\mathcal{L}\otimes I_n)(x(t)+e_d(t))$. Choose, without loss of generality, agent $r=1$ and define $\eta(t)=x_i(t)-x_1(t)$. It follows that
\begin{align}
\left.
	\begin{array}{l l}
 \dot{\eta}(t)\!&=\! -(\mathcal{L}_s\otimes I_n)\eta(t)-(\mathcal{L}'\otimes I_n)e_d(t).
			\end{array}  \label{eq:Tranfdynco}  \right.
\end{align}
Then, we have that 
\begin{align}
  \left.
	\begin{array}{l l}
	\left\|\eta(t)\right\| \leq \hat{\beta}\eta_0\textbf{e}^{-\hat{\lambda} t} + \frac{\hat{\beta}\gamma\sqrt{N}L'}{\hat{\lambda}-\lambda} (\textbf{e}^{-\lambda t}-\textbf{e}^{-\hat{\lambda} t}) 
\end{array}   \label{eq:x2hatresco} \right.
\end{align}
and the agents achieve consensus asymptotically. In this case we can see that the initial average remains constant
\begin{align}
  \left.
	\begin{array}{l l}
	\frac{d}{dt}\frac{1}{N}1_N^Tx(t)=-\frac{1}{N}1_N^T (\mathcal{L}\otimes I_n)(x+\xi)=0
	\end{array}    \right.
\end{align}
since the graph is undirected; thus, the group of agents reach  average consensus. 

Finally, we use the expression in the right-hand side of \eqref{eq:x2hatresco} in order to bound the growth of the error $e_i$ between triggering events as it was done in \eqref{eq:errdyn}-\eqref{eq:errresp4-1} to obtain the lower-bound \eqref{eq:tauco} on the inter-event time intervals. The admissible delay \eqref{eq:delayco} is obtained in a similar way.

 $\square$

\section{Example} \label{sec:example}
Let us consider six agents connected using a directed graph. The adjacency matrix $\mathcal{A}$ is given by $a_{12}=a_{23}=a_{25}=a_{32}=a_{36}=a_{43}=a_{45}=a_{54}=a_{61}=1$ and the remaining entries of $\mathcal{A}$ are equal to zero.
The communication channel is subject to packet dropouts and the MANSD is $\rho-1$, where $\rho=4$. The initial conditions are $x(0)=[1\ -\!1\ 2\ 3\ 5\ 4]^\text{T}$. We choose the parameters $\beta=1$, $\lambda=0.4$, $\gamma_d=9$, and $\bar{\delta}=1.5$. The obtained minimum inter event time is $\tau=0.0025$ and the admissible delays are bounded by $d=0.0223$. The values of $\tau$ and $d$ represent stationary values, when $t_{k_i}\rightarrow \infty$, and they also represent the minimum values for any $t_{k_i}\in[0,\infty)$. Note that in the simulation, the error events are not continuously evaluated but a sampling time $\tau_s<\tau$ is implemented, where $\tau_s=0.0002$. Future work will consider the effects of non-continuous event detection.

\begin{figure}
	\begin{center}
		\includegraphics[width=12cm,height=6.2cm,trim=.9cm .1cm .9cm .4cm]{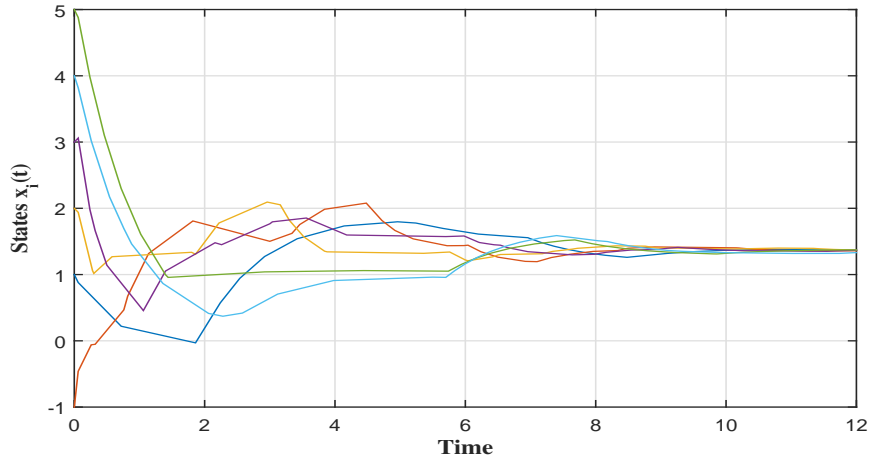}
	\caption{States of the $6$ agents}
	\label{fig:states}
	\end{center}
\end{figure}

\begin{figure}
	\begin{center}
		\includegraphics[width=12cm,height=7.5cm,trim=.9cm .6cm .9cm .6cm]{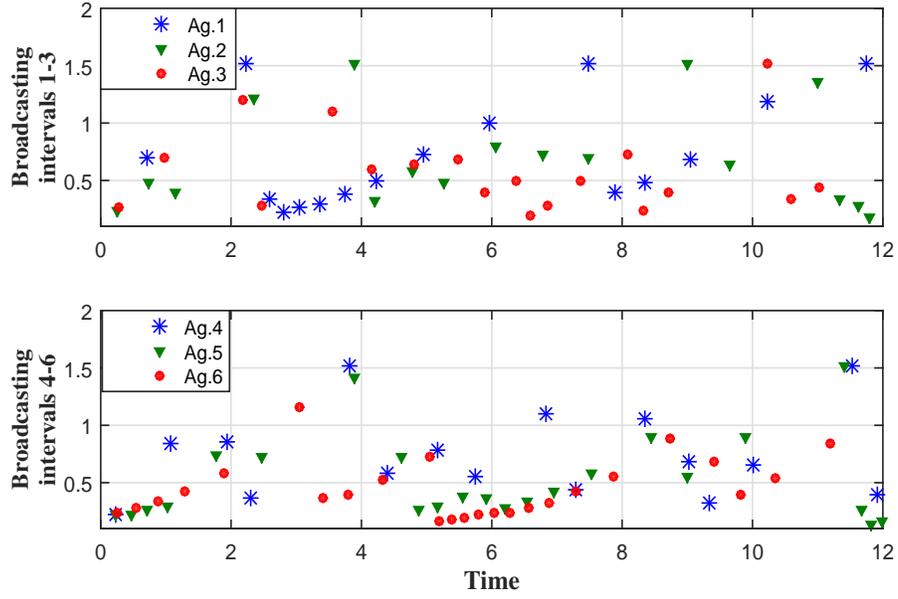}
	\caption{Broadcasting time intervals for each agent}
	\label{fig:Bro}
	\end{center}
\end{figure}

\begin{figure}
	\begin{center}
		\includegraphics[width=12cm,height=7.5cm,trim=.9cm .6cm .9cm .1cm]{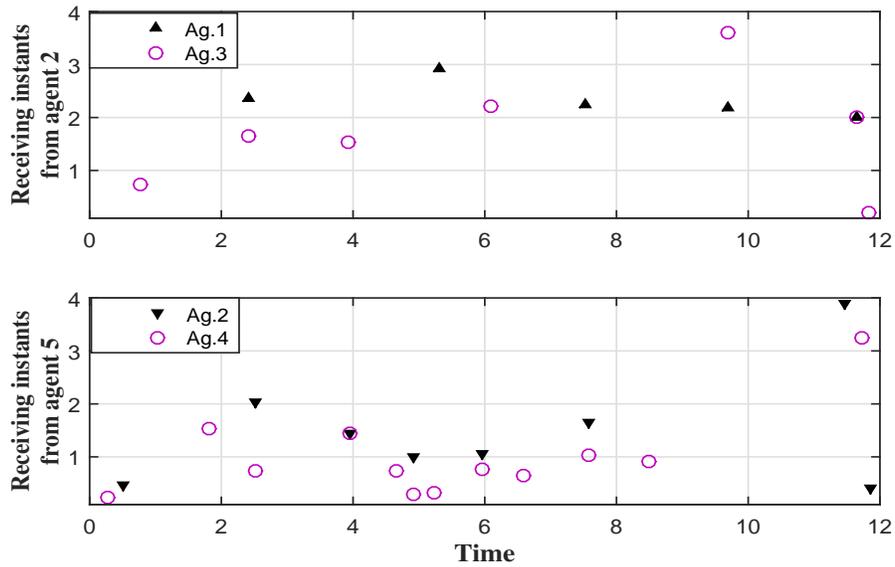}
	\caption{Receiving time intervals for agents $j$ such that $2\in\mathcal{N}_j$ (top) and $5\in\mathcal{N}_j$ (bottom)}
	\label{fig:Rec}
	\end{center}
\end{figure}

Fig. \ref{fig:states} shows the response of the six agents for time-varying delays $d_{ij}(t_{k_i})\in[0.005,0.02]$. Additionally, measurement updates generated by agent $i$ based on its local events may be lost and may not be received by some or by all of the intended agents $j$, $i\in\mathcal{N}_j$, and for $i=1,...,N$. The maximum number of successive dropped packets is $\rho-1$. It can be seen from Fig. \ref{fig:states} that all agents converge to a common value, that is, consensus is reached. The time intervals between events are shown for each one of the agents in Fig. \ref{fig:Bro}; however, some of these updates do not reach their destinations. Fig. \ref{fig:Rec} shows the receiving time intervals from agent 2 to agents 1 and 3; it also shows the receiving time intervals from agent 5 to agents 2 and 4. It can be clearly seen, for instance, that only a fraction of the measurement updates generated by agent 2 are able to reach agents 1 or 3. Thus, the corresponding receiving time intervals are much greater, in general, than the broadcasting time intervals. Similar situations occur to every agent receiving information with packet losses.

\section{Conclusions} \label{sec:conclusion}
This paper presented a decentralized event-triggered consensus protocol for a group of agents connected using directed graphs where the communication channel is subject to non-consistent packet dropouts and non-consistent communication delays. It was shown that agents are able to achieve asymptotic consensus under this protocol and that Zeno behavior is avoided. The use of event-triggered control techniques is beneficial when continuous communication among agents is not possible and it provides a higher level of decentralization since it is not necessary for agents to know a global sampling period and global communication time instants as in sampled-data approaches. The decentralized event-triggered consensus protocol provides each agent the freedom to decide their own broadcasting instants independently of any other agent in the network. Future work will consider the cases of double integrator dynamics and agents with linear dynamics under packet loss and communication delays.


\begin{thebibliography}{10}

\bibitem{Ji07}
M.~Ji and M.~Egerstedt, ``Distributed coordination control of multiagent
  systems while preserving connectedness,'' \emph{IEEE Transactions on
  Robotics}, vol.~23, no.~4, pp. 693--703, 2007.

\bibitem{Moreau04}
L.~Moreau, ``Stability of continuous-time distributed consensus algorithms,''
  in \emph{IEEE Conference on Decision and Control}, 2004, pp. 3998--4003.

\bibitem{Olfati04}
R.~Olfati-Saber and R.~M. Murray, ``Consensus problems in networks of agents
  with switching topology and time-delays,'' \emph{IEEE Transactions on
  Automatic Control}, vol.~49, no.~9, pp. 1520--1533, 2004.

\bibitem{CaoRen10}
Y.~Cao and W.~Ren, ``Multi-vehicle coordination for double integrator dynamics
  under fixed undirected/directed interaction in a sampled-data setting,''
  \emph{International Journal of Robust and Nonlinear Control}, vol.~20, no.~9,
  pp. 987--1000, 2010.

\bibitem{QinGao12}
J.~Qin and H.~Gao, ``A sufficient condition for convergence of sampled-data
  consensus for double integrator dynamics with nonuniform and time-varying
  communication delays,'' \emph{IEEE Transactions on Automatic Control},
  vol.~57, no.~9, pp. 2417--2422, 2012.

\bibitem{Hayakawa06}
T.~Hayakawa, T.~Matsuzawa, and S.~Hara, ``Formation control of multi-agent
  systems with sampled information,'' in \emph{45th IEEE Conference on Decision
  and Control}, 2006, pp. 4333--4338.

\bibitem{Liu10}
H.~Liu, G.~Xie, and L.~Wang, ``Necessary and sufficient conditions for solving
  consensus of double integrator dynamics via sampled control,''
  \emph{International Journal of Robust and Nonlinear Control}, vol.~20,
  no.~15, pp. 1706--1722, 2010.

\bibitem{Garcia13b}
E.~Garcia, Y.~Cao, H.~Yu, P.~J. Antsaklis, and D.~W. Casbeer, ``Decentralized
  event-triggered cooperative control with limited communication,''
  \emph{International Journal of Control}, vol.~86, no.~9, pp. 1479--1488,
  2013.

\bibitem{kia14}
S.~Kia, J.~Cort{\'e}s, and S.~Mart{\'\i}nez, ``Distributed event-triggered
  communication for dynamic average consensus in networked systems,''
  \emph{Automatica}, vol.~59, pp. 112--119, 2015.

\bibitem{Meng13}
X.~Meng and T.~Chen, ``Event based agreement protocols for multi-agent
  networks,'' \emph{Automatica}, vol.~49, no.~7, pp. 2125--2132, 2013.

\bibitem{Nowzari14}
C.~Nowzari and J.~Cortes, ``Zeno-free, distributed event-triggered
  communication and control for multi-agent average consensus,'' in
  \emph{American Control Conference}, 2014, pp. 2148--2153.

\bibitem{Seyboth13}
G.~S. Seyboth, D.~V. Dimarogonas, and K.~H. Johansson, ``Event-based
  broadcasting for multi-agent average consensus,'' \emph{Automatica}, vol.~49,
  no.~1, pp. 245--252, 2013.

\bibitem{AntaTabuada10}
A.~Anta and P.~Tabuada, ``To sample or not to sample: Self-triggered control
  for nonlinear systems,'' \emph{IEEE Transactions on Automatic Control},
  vol.~55, no.~9, pp. 2030--2042, 2010.

\bibitem{Astrom02}
K.~J. Astrom and B.~M. Bernhardson, ``Comparison of {R}iemann and {L}ebesgue
  sampling for first order stochastic systems,'' in \emph{41st IEEE Conference
  on Decision and Control}, 2002, pp. 2011--2016.

\bibitem{Garcia13}
E.~Garcia and P.~J. Antsaklis, ``Model-based event-triggered control for
  systems with quantization and time-varying network delays,'' \emph{IEEE
  Transactions on Automatic Control}, vol.~58, no.~2, pp. 422--434, 2013.
	
\bibitem{Donkers12}
M.~Donkers and W.~Heemels, ``Output-based event-triggered control with
  guaranteed-gain and improved and decentralized event-triggering,''
  \emph{Automatic Control, IEEE Transactions on}, vol.~57, no.~6, pp.
  1362--1376, 2012.

\bibitem{Garcia12IJC}
E.~Garcia and P.~J. Antsaklis, ``Parameter estimation and adaptive stabilization in
time-triggered and event-triggered model-based control of uncertain systems,'' \emph{
International Journal of Control}, vol.~85, no.~9, pp. 1327--1342, 2012.

\bibitem{Tabuada07}
P.~Tabuada, ``Event-triggered real-time scheduling of stabilizing control
  tasks,'' \emph{IEEE Transactions on Automatic Control}, vol.~52, no.~9, pp.
  1680--1685, 2007.

\bibitem{WangLemmon11}
X.~Wang and M.~Lemmon, ``Event-triggering in distributed networked control
  systems,'' \emph{IEEE Transactions on Automatic Control}, vol.~56, no.~3, pp.
  586--601, 2011.
	
\bibitem{Garcia12CDC}
E.~Garcia and P.~J. Antsaklis, ``Output feedback model-based control of uncertain
discrete-time systems with network induced delays,'' in \emph{51st
IEEE Conference on Decision and Control}, 2012, pp. 6647--6652.

\bibitem{Garcia14CDCdis}
E.~Garcia and P.~J. Antsaklis, ``Event-triggered output feedback stabilization of networked 
systems with external disturbance,'' in \emph{53st
IEEE Conference on Decision and Control}, 2014, pp. 3566--3671.

\bibitem{Guinaldo13}
M.~Guinaldo, D.~Dimarogonas, K.~H. Johansson, J.~S{\'a}nchez, and S.~Dormido,
  ``Distributed event-based control strategies for interconnected linear
  systems,'' \emph{Control Theory \& Applications, IET}, vol.~7, no.~6, pp.
  877--886, 2013.

\bibitem{Guinaldo14}
M.~Guinaldo, D.~Lehmann, J.~S{\'a}nchez, S.~Dormido, and K.~H. Johansson,
  ``Distributed event-triggered control for non-reliable networks,''
  \emph{Journal of the Franklin Institute}, vol. 351, no.~12, pp. 5250--5273,
  2014.

\bibitem{Mazo11TAC}
M.~Mazo and P.~Tabuada, ``Decentralized event-triggered control over wireless
  sensor/actuator networks,'' \emph{IEEE Transactions on Automatic Control},
  vol.~56, no.~10, pp. 2456--2461, 2011.

\bibitem{wang2010relaxing}
X.~Wang, Y.~Sun, and N.~Hovakimyan, ``Relaxing the consistency condition in
  distributed event-triggered networked control systems,'' in \emph{49th IEEE
  Conference on Decision and Control}, 2010, pp. 4727--4732.

\bibitem{ChenHao12}
X.~Chen and F.~Hao, ``Event-triggered average consensus control for
  discrete-time multi-agent systems,'' \emph{IET Control Theory and
  Applications}, vol.~6, no.~16, pp. 2493--2498, 2012.

\bibitem{YinYue13b}
X.~Yin and D.~Yue, ``Event-triggered tracking control for heterogeneous
  multi-agent systems with {M}arkov comunication delays,'' \emph{Journal of the
  Franklin Institute}, vol. 350, no.~5, pp. 1312--1334, 2013.

\bibitem{GuoDimarogonas13}
M.~Guo and D.~V. Dimarogonas, ``Nonlinear consensus via continuous, sampled,
  and aperiodic updates,'' \emph{International Journal of Control}, vol.~86,
  no.~4, pp. 567--578, 2013.

\bibitem{Garcia14Auto}
{E. Garcia, Y. Cao and D. W. Casbeer}, ``Decentralized event-triggered
  consensus with general linear dynamics,'' \emph{Automatica}, vol.~50, no.~10,
  pp. 2633--2640, 2014.

\bibitem{Garcia14}
E.~Garcia, Y.~Cao, and D.~W. Casbeer, ``Cooperative control with general linear
  dynamics and limited communication: centralized and decentralized
  event-triggered control strategies,'' in \emph{American Control Conference},
  2014, pp. 159--164.

\bibitem{LiuHill12}
T.~Liu, D.~J. Hill, and B.~Liu, ``Synchronization of dynamical networks with
  distributed event-based communication,'' in \emph{51st IEEE Conference on
  Decision and Control}, 2012, pp. 7199--7204.

\bibitem{Zhu14}
W.~Zhu, Z.-P. Jiang, and G.~Feng, ``Event-based consensus of multi-agent
  systems with general linear models,'' \emph{Automatica}, vol.~50, no.~2, pp.
  552--558, 2014.

\bibitem{Garcia14CDC}
E.~Garcia, Y.~Cao, and D.~W. Casbeer, ``Event-triggered cooperative control
  with general linear dynamics and communication delays,'' in \emph{IEEE
  Conference on Decision and Control}, 2014, pp. 2914--2919.

\bibitem{kia15}
S.~Kia, J.~Cort{\'e}s, and S.~Mart{\'\i}nez, ``Distributed convex optimization
  via continuous-time coordination algorithms with discrete-time
  communication,'' \emph{Automatica}, vol.~55, pp. 254--264, 2015.

\bibitem{Aragues12}
R.~Aragues, G.~Shi, D.~V. Dimarogonas, C.~Sagues, and K.~H. Johansson,
  ``Distributed algebraic connectivity estimation for adaptive event-triggered
  consensus,'' in \emph{American Control Conference}, 2012, pp. 32--37.

\bibitem{Franceschelli09}
M.~Franceschelli, A.~Gasparri, A.~Giua, and C.~Seatzu, ``Decentralized
  laplacian eigenvalues estimation for networked multi-agent systems,'' in
  \emph{48th IEEE Conference on Decision and Control and 28th Chinese Control
  Conference}, 2009, pp. 2717--2722.

\end{thebibliography}
\end{document}